\documentclass[preprint,12pt]{elsarticle}

\usepackage{hyperref}
\usepackage{graphicx}

\journal{Journal of Applied Numerical Mathematics}

\begin{document}

\begin{frontmatter}

\title{On the kernel and particle consistency in smoothed particle hydrodynamics}


\author[address1,address2]{Leonardo Di G. Sigalotti\corref{mycorrespondingauthor}}
\cortext[mycorrespondingauthor]{Corresponding author}
\ead{leonardo.sigalotti@gmail.com}

\author[address3,address4]{Jaime Klapp}
\ead{jaime.klapp@inin.gob.mx}

\author[address2]{Otto Rend\'on}
\ead{ottorendon@gmail.com}

\author[address1]{Carlos A. Vargas}
\ead{carlosvax@gmail.com}

\author[address4,address2]{Franklin Pe\~na-Polo}
\ead{franklin.pena@gmail.com}

\address[address1]{\'Area de F\'{\i}sica de Procesos Irreversibles, Departamento de
Ciencias B\'asicas, Universidad Aut\'onoma Metropolitana-Azcapotzalco (UAM-A),
Av. San Pablo 180, 02200 M\'exico D.F., Mexico}
\address[address2]{Centro de F\'{\i}sica, Instituto Venezolano de Investigaciones
Cient\'{\i}ficas (IVIC), Apartado Postal 20632, Caracas 1020-A, Venezuela}
\address[address3]{Departamento de F\'{\i}sica, Instituto Nacional de Investigaciones
Nucleares (ININ), Carretera M\'exico-Toluca km. 36.5, La Marquesa, 52750 Ocoyoacac, Estado de 
M\'exico, Mexico}
\address[address4]{ABACUS-Centro de Matem\'aticas Aplicadas y C\'omputo de Alto
Rendimiento, Departamento de Matem\'aticas, Centro de Investigación y de Estudios
Avanzados (Cinvestav-IPN), Carretera M\'exico-Toluca km. 38.5, La Marquesa, 52740
Ocoyoacac, Estado de M\'exico, Mexico}

\begin{abstract}
The problem of consistency of smoothed particle hydrodynamics (SPH) has demanded 
considerable attention in the past few years due to the ever increasing number of
applications of the method in many areas of science and engineering. A loss of
consistency leads to an inevitable loss of approximation accuracy. In this paper,
we revisit the issue of SPH kernel and particle consistency and demonstrate that
SPH has a limiting second-order convergence rate. Numerical experiments with suitably
chosen test functions validate this conclusion. In particular, we find that when using
the root mean square error as a model evaluation statistics, well-known corrective
SPH schemes, which were thought to converge to second, or even higher order,
are actually first-order accurate, or at best close to second order. We also
find that observing the joint limit when $N\to\infty$, $h\to 0$, and $n\to\infty$,
as was recently proposed by Zhu et al., where $N$ is the total number of particles,
$h$ is the smoothing length, and $n$ is the number of neighbor particles, standard
SPH restores full $C^{0}$ particle consistency for both the estimates of the function
and its derivatives and becomes insensitive to particle disorder. 
\end{abstract}


\end{frontmatter}


\section{Introduction}

The method of smoothed particle hydrodynamics (SPH) was introduced in the literature
independently by Gingold and Monaghan \cite{Gingold77} and Lucy \cite{Lucy77} for modeling
astrophysical flow problems. Since then the method has been widely applied to different
areas of science and engineering due to its simplicity and ease of implementation. At
the same time, the method has been improved over the years to overcome major shortcomings
and deficiencies. One longstanding drawback of standard SPH is the particle inconsistency, 
which is an intrinsic manifestation of the lack of integrability of the kernel
approximation in its spatially discretized form, resulting in a loss of accuracy.
In practical applications, SPH inconsistency arises when the support domain of the kernel
is truncated by a model boundary, for irregularly distributed particles even in the
absence of kernel truncation, and for spatially adaptive calculations where a variable
smoothing length is employed.

Several corrective strategies have been proposed to restore particle consistency in SPH
calculations. A simple correction technique was first advanced by Li and Liu \cite{Li96}
and Liu at al. \cite{Liu97}, where the kernel itself is modified to ensure that polynomial
functions up to a given degree are exactly interpolated. A kernel gradient correction,
allowing for the exact evaluation of the gradient of a linear function, was further
proposed by Bonet and Lok \cite{Bonet99} based on a variational formulation of SPH. A
general approach to the construction of kernel functions that obey the consistency
conditions of SPH in continuous form and describe the compact supportness requirement
was presented by Liu et al. \cite{Liu03}. A drawback of this approach is that the
reconstructed smoothing functions may be partially negative, non-symmetric, and 
non-monotonically decreasing, thereby compromising the stability of the numerical
simulations. 

More stable approaches for restoring consistency are based on Taylor
series expansions of the kernel approximations of a function and its derivatives. In
general, if up to $m$ derivatives are retained in the series expansions, the resulting
kernel and particle approximations will have $(m+1)$th-order accuracy or $m$th-order
consistency (i.e., $C^{m}$ consistency). This approach was first developed by
Chen et al. \cite{Chen99a,Chen99b} (their corrective smoothed particle method, or
CSPM), which solves for the approximation of a function separately from that of its 
derivatives by neglecting all terms involving derivatives in the former expansion and 
retaining only first-order terms in the latter expansions. This scheme is equivalent
to a Shepard's interpolation for the function \cite{Randles96}, and so it should
restore $C^{1}$ kernel and particle consistency for the interior regions and $C^{0}$ 
consistency at the boundaries. By retaining only first-order derivatives in the Taylor
series expansions for the function and its derivatives and solving simultaneously the
resulting set of linear equations (the FPM scheme), Liu and Liu \cite{Liu06} argued
that $C^{1}$ kernel and particle consistency can be obtained for both interior and
boundary regions. The FPM scheme was further improved by Zhang and Batra \cite{Zhang04} 
(their MSPH scheme), where now up to second-order derivatives are retained in the
Taylor series expansions. In principle, this method should restore $C^{2}$ consistency
for the SPH approximation of the function (i.e., third-order convergence rates) and
$C^{1}$ consistency for the first-order derivatives. However, when adding higher-order
derivatives in the Taylor expansions, the number of algebraic linear equations to be
solved increases rapidly, implying high computational costs. In addition, since the
solution involves a matrix inversion, for some types of problems the stability of the
scheme can be compromised by the conditioning of the matrix. A modified FPM approach,
which is free of kernel gradients and leads to a symmetric corrective matrix was 
recently proposed by Huang et al. \cite{Huang15}. An alternative formulation based on
the inclusion of boundary integrals in the kernel approximation of the spatial 
derivatives was reported by Maci\`{a} et al. \cite{Macia12}, which restores $C^{0}$
consistency at the model boundaries. A new SPH formulation, based on a novel piecewise
high-order Moving-Least-Squares WENO reconstruction and on the use of Riemann solvers,
that improves the accuracy of SPH in the presence of sharp discontinuities was recently
reported by Avesani et al. \cite{Avesani14}.  

A new strategy to ensure formal convergence and particle consistency with standard
SPH has recently been devised by Zhu et al. \cite{Zhu15} in the astrophysical
context. In this approach, no corrections are required and full consistency is
recovered provided that $N\to\infty$, $h\to 0$, and $n\to\infty$, where $N$ is
the total number of particles, $h$ is the smoothing length, and $n$ is the number
of neighbor particles within the kernel support. They found that if $n$ is
fixed, as is customary in SPH, consistency will not be guaranteed even though
$N\to\infty$ and $h\to 0$ since there is a residual error that will not vanish
unless $n$ is allowed to increase with $N$ as $n\sim N^{1/2}$. However, the systematic 
increase of $n$ with improved resolution demands changing the interpolation
kernel to a Wendland-type function \cite{Wendland95}, which, unlike traditional
kernels, is free from the pairing instability when used to perform smoothing
in SPH with large numbers of neighbors \cite{Dehnen12}.

In this paper, we revisit the issue of kernel and particle consistency in SPH. We
first demonstrate that the normalization condition of the kernel is independent of
$h$, suggesting that its discrete representation depends only on $n$, consistently
with the error analyses of Vaughan et al. \cite{Vaughan08} and Read et al. \cite{Read10}. 
Although $C^{0}$ and $C^{1}$ kernel and particle consistency can be achieved by
some corrective SPH methods, a simple observation shows that $C^{2}$ kernel consistency
is difficult to achieve, implying an upper limit to the convergence rate of SPH in
practical applications. Numerical experiments with suitably chosen test functions 
in two-space dimensions validate this conclusion. The paper is organized as follows.
In Section 2, we discuss the issue of $C^{0}$ consistency and show that the
normalization condition of the kernel is independent of $h$. The issue of higher-order
consistency is considered in Section 3, where we show that $C^{2}$ consistency is
affected by an inherent intrinsic diffusion, which arises as a consequence of the
dispersion of the SPH particle positions relative to the mean. Section 4 outlines the 
importance of restoring $C^{1}$ consistency for the gradient. Finally, Section 5 presents 
numerical tests that demonstrate the convergence rates of the particle approximations and
Section 6 contains the conclusions.

\section{Normalization condition and $C^{0}$ consistency}

As it is well-known, the starting point of SPH lies on the exact identity
\begin{equation}
f({\bf x})=\int _{{\cal R}^{3}}f({\bf x}^{\prime})\delta ({\bf x}-{\bf x}^{\prime})
d{\bf x}^{\prime},
\end{equation}
where $f=f({\bf x})$ is some sufficiently smooth function,
$\delta ({\bf x}-{\bf x}^{\prime})$ is the Dirac-$\delta$ distribution, and the
integration is taken over all space. The kernel approximation is obtained by
replacing the Dirac-$\delta$ distribution by some kernel interpolation function $W$
such that
\begin{equation}
\langle f({\bf x})\rangle =\int _{{\cal R}^{3}}f({\bf x}^{\prime})
W(|{\bf x}-{\bf x}^{\prime}|,h)d{\bf x}^{\prime},
\end{equation}
where $\langle f({\bf x})\rangle$ is the smoothed estimate of $f({\bf x})$ and 
$W$ must satisfy the following properties: (a) in the limit $h\to 0$ it becomes
the Dirac-$\delta$ function so that $\langle f({\bf x})\rangle\to f({\bf x})$,
(b) it must satisfy the normalization condition
\begin{equation}
\int _{{\cal R}^{3}}W(|{\bf x}-{\bf x}^{\prime}|,h)d{\bf x}^{\prime}=1,
\end{equation}
and (c) it must be positive definite, symmetric, and monotonically decreasing. Almost
all modern applications of SPH assume that $W$ has compact support, that is,
$W(|{\bf x}-{\bf x}^{\prime}|,h)=0$ for $|{\bf x}-{\bf x}^{\prime}|\geq kh$,
where $k$ is some integer that depends on the kernel function.

We shall first demonstrate that the normalization condition for a kernel
satisfying the above properties is independent of the smoothing length $h$. This
feature is tacitly assumed in the SPH literature. However, as we shall
see later, it has important conceptual implications for the kernel consistency
relations. For simplicity in exposition, let us restrict to one-space dimension and
assume the kernel to be a Gaussian function such that
\begin{equation}
\delta (x-x^{\prime})=\lim _{h\to 0}\frac{1}{\sqrt{2\pi}h}\exp\left[
\frac{-\left(x-x^{\prime}\right)^{2}}{2h^{2}}\right].
\end{equation}
Making the change $|x-x^{\prime}|\to h|x-x^{\prime}|$ in the Gaussian kernel, it
is then easy to show the following scaling relation
\begin{equation}
W(h|x-x^{\prime}|,h)=\frac{1}{h}W(|x-x^{\prime}|,1),
\end{equation}
which is indeed satisfied by all known SPH kernel functions. Similar forms follow
in two- and three-space dimensions with $1/h$ replaced by $1/h^{2}$ and $1/h^{3}$,
respectively. Now expanding in Taylor series $f(x^{\prime})$ around $x^{\prime}=x$
with $x\to hx$ and $x^{\prime}\to hx^{\prime}$, we obtain
\begin{equation}
f(hx^{\prime})=\sum _{l=0}^{\infty}\frac{1}{l!}h^{l}f^{(l)}(hx)(x^{\prime}-x)^{l},
\end{equation}
which when inserted in the kernel approximation (2) yields
\begin{equation}
\langle f(hx)\rangle =\sum _{l=0}^{\infty}\frac{1}{l!}h^{l}f^{(l)}(hx)
\int _{\cal R}(x^{\prime}-x)^{l}W(|x-x^{\prime}|,1)dx^{\prime},
\end{equation}
where we have made $dx^{\prime}\to hdx^{\prime}$ in the integrand and used the
scaling relation (5). Using the Gaussian kernel, the integral in Eq. (7) when
$l=0$ is just the probability function and is exactly one, yielding the
condition
\begin{equation}
\int _{\cal R}W(|x-x^{\prime}|,1)dx^{\prime}=1.
\end{equation}
Therefore, the normalization condition is independent of $h$ provided that the
kernel interpolation obeys the scaling relation (5). Alternatively, in the SPH
literature the kernels are usually defined in terms of the dimensionless parameter
$q=|x-x^{\prime}|/h$ such that $W(|x-x^{\prime}|,h)=W(q)/h$, which is not the same
as Eq. (5) because $q$ depends explicitly on $h$. In addition, all terms in Eq. (7)
with $l$ odd vanish because of the symmetry of the kernel function, while only those
with $l$ even survive. Therefore, using Eq. (8) and retaining only the $l=2$ term in
Eq. (7), we obtain
\begin{equation}
\langle f(hx)\rangle =f(hx)+\frac{1}{2}h^{2}f^{\prime\prime}(hx)\int _{\cal R}
(x^{\prime}-x)^{2}W(|x-x^{\prime}|,1)dx^{\prime}+O(h^{4}),
\end{equation}
which expresses that the kernel approximation of a function is second-order accurate
and therefore has $C^{1}$ consistency for an unbounded domain. Thus, for such infinitely
extended domains $C^{1}$ consistency requires that $C^{0}$ consistency be satisfied
according to Eq. (8). As was pointed out by Liu and Liu \cite{Liu06} and formally 
demonstrated by Vaughan \cite{Vaughan09}, for a bounded domain the kernel
approximation (2) needs be replaced by the integral form of the Shepard interpolant
to guarantee $C^{0}$ consistency near a model boundary. Equation (9) bears some
resemblance with the expression of the error derived by Vaughan et al. \cite{Vaughan08},
where the contribution to the error due to the smoothing length can be separated from
that due to the discretized form of the integral, which, being independent of $h$, will
directly depend on the number of neighbors within the kernel support. 

The SPH discretization makes reference to a set of Lagrangian particles which may, in
general, be disordered. If we consider a finite model domain $\Omega\in {\cal R}$ and
divide it into $N$ sub-domains, labeled $\Omega _{a}$, each of which contains a 
Lagrangian particle $a$ at position $x_{a}\in\Omega _{a}$, the discrete form of 
Eq. (8) becomes
\begin{equation}
\sum _{b=1}^{N}W_{ab}\Delta x_{b}=O(1),
\end{equation}
where $W_{ab}=W(|x_{a}-x_{b}|,1)$. In general, the above summation is not exactly one.
In this approximation, the error depends on the number of particles $N$ and how these
are actually distributed. If the kernel has compact support, then the above summation
is over the number of particles $n$ within a length over which the kernel itself does
not vanish, i.e., $|x_{a}-x_{b}|<kh$. In this case, the error scales as $\sim n^{-\alpha}$,
where $\alpha =0.5$ for a truly random distribution of particles and $\alpha =1$
for quasi-ordered patterns \cite{Zhu15}. In either case, if $n$ is sufficiently
large the discrete normalization condition becomes sufficiently close to one ensuring
$C^{0}$ particle consistency. This will certainly require to scale $n$ with $N$ as 
$\sim N^{0.5}$ as the resolution is improved \cite{Zhu15}.

The discrete form (10) losses $C^{0}$ consistency for irregularly distributed particles
even far away from the model boundaries because node disorder results in a noise
error, which scales with $n$ as $\sim n^{-1}\ln n$ \cite{Monaghan92}. $C^{0}$ particle
consistency is equivalent to demand that the homogeneity of space is not affected by
the process of spatial discretization, which in turn has, as a consequence, the
conservation of linear momentum \cite{Vignjevic09}. In other words, the SPH interpolation
has to be independent of a rigid-body translation of the coordinates. To see this,
consider the discrete form of Eq. (2) for the position vector ${\bf x}=(x,y,z)$, i.e.,
\begin{equation}
\langle {\bf x}\rangle _{a}=\sum _{b=1}^{n}{\bf x}_{b}W_{ab}\Delta V_{b},
\end{equation}
where $W_{ab}=W(|{\bf x}_{a}-{\bf x}_{b}|,1)$ and $\Delta V_{b}$ is the volume of the
sub-domain $\Omega _{b}$ within which particle $b$ lies. In SPH simulations, it is
common practice to evaluate $\Delta V_{b}$ as the ratio $m_{b}/\rho _{b}$, where
$m_{b}$ and $\rho_{b}$ refer to the mass and density of particle $b$, respectively.
According to Eq. (11), the estimate of the transformed coordinates 
${\bf x}^{\prime}={\bf x}+\Delta {\bf x}$ is
\begin{equation}
\langle {\bf x}^{\prime}\rangle _{a}=\sum _{b=1}^{n^{\prime}}{\bf x}_{b}^{\prime}
W_{ab}^{\prime}\Delta V_{b}^{\prime},
\end{equation}
where $W_{ab}^{\prime}=W(|{\bf x}_{a}^{\prime}-{\bf x}_{b}^{\prime}|,1)$.
Preservation of space homogeneity under uniform translation demands that
$\langle\Delta {\bf x}\rangle =\Delta {\bf x}$ so that
$\langle {\bf x}^{\prime}\rangle _{a}=\langle {\bf x}\rangle _{a}+\Delta {\bf x}$.
Replacing ${\bf x}_{b}^{\prime}$ by ${\bf x}_{b}+\Delta {\bf x}$ in Eq. (12) yields
\begin{equation}
\langle {\bf x}^{\prime}\rangle _{a}=\langle {\bf x}\rangle _{a}+
\Delta {\bf x}\sum _{b=1}^{n}W_{ab}\Delta V_{b},
\end{equation}
where we have made $W_{ab}^{\prime}\Delta V_{b}^{\prime}=W_{ab}\Delta V_{b}$
and $n^{\prime}=n$ because under solid-body translation the coordinates of a point
are independent of the translation of the coordinates axes. Therefore, Eq. (13)
expresses that homogeneity of the discretized space is recovered by the SPH 
interpolation only if the condition
\begin{equation}
\sum _{b=1}^{n}W_{ab}\Delta V_{b}=1,
\end{equation}
is satisfied exactly. A simple way to enforce this condition at a model boundary is 
to use a discrete Shepard interpolation \cite{Randles96,Bonet99,Vaughan09}, where 
the kernel function is normalized according to
\begin{equation}
W_{ab}\to\frac{W_{ab}}{\sum _{b=1}^{n}W_{ab}\Delta V_{b}},
\end{equation}
or alternatively, to use huge numbers of particles such that the limits $N\to\infty$, 
$h\to 0$, and $n\to\infty$ are achieved in an approximate sense \cite{Zhu15}.

\section{Higher-order consistency}

For the sake of simplicity we shall restrict ourselves again to one-space dimension.
However, the analysis of this and the preceding section can be easily generalized to
higher dimensions. If a polynomial of degree $m$ is exactly reproduced by an
approximation, then we say that it has $C^{m}$ consistency and the following family
of integral conditions must hold for the kernel
\begin{equation}
\int _{\cal R}(x-x^{\prime})^{l}W(|x-x^{\prime}|,h)dx^{\prime}=\delta_ {l0},
\end{equation}
for $l=0$,1,2,...$m$, where $\delta _{l0}$($=1$ for $l=0$ and $=0$ for $l\neq 0$) is
the Kronecker delta. Using the scaling relation (5) and making 
$|x-x^{\prime}|\to h|x-x^{\prime}|$, the integrals in Eq. (16) can be rewritten as
\begin{equation}
\int _{\cal R}(x-x^{\prime})^{l}W(|x-x^{\prime}|,1)dx^{\prime}=h^{-l}\delta_ {l0}.
\end{equation}
For $l=0$, Eq. (17) states the normalization condition (8), while for $l=1$ it
states the symmetry property of the kernel. Using the Gaussian kernel, we may show
that this integral is exactly zero for an infinitely extended domain, implying 
$C^{1}$ kernel consistency. Thus, if $C^{1}$ consistency is achieved, then $C^{0}$
consistency is automatically guaranteed.

For $l=2$, Eq. (17) reduces to
\begin{equation}
\int _{\cal R}(x-x^{\prime})^{2}W(|x-x^{\prime}|,1)dx^{\prime}=0,
\end{equation}
which must be satisfied in order to ensure $C^{2}$ consistency for the kernel
approximation. We note that this integral already appears in Eq. (9) as a finite
source of error and it does not vanish unless the kernel approaches the Dirac-$\delta$
function. Using Eqs. (2) and (8), we may evaluate this integral to obtain
\begin{eqnarray}
\int _{\cal R}(x-x^{\prime})^{2}W(|x-x^{\prime}|,1)dx^{\prime}&=&
x^{2}-2x\langle x\rangle +\langle x^{2}\rangle\nonumber\\
&=&\left(x-\langle x\rangle\right)^{2}+\left(\langle x^{2}\rangle -
\langle x\rangle ^{2}\right).
\end{eqnarray}
If $C^{1}$ consistency is guaranteed by the kernel approximation then $x=\langle x\rangle$,
and the above integral becomes
\begin{equation}
\int _{\cal R}(x-x^{\prime})^{2}W(|x-x^{\prime}|,1)dx^{\prime}=
\langle x^{2}\rangle -\langle x\rangle ^{2}\neq 0,
\end{equation}
which implies that $C^{2}$ kernel consistency is not achieved even though $C^{0}$ and
$C^{1}$ consistencies are satisfied, unless $W(|x-x^{\prime}|,1)\to\delta (x-x^{\prime})$,
or equivalently, $N\to\infty$, $h\to 0$, and $n\to\infty$ as suggested by Zhu et al.
\cite{Zhu15}, in which case $\langle x^{2}\rangle =\langle x\rangle ^{2}$.
It is evident from Eq. (20) that the lack of $C^{2}$ consistency is due to an 
intrinsic diffusion $\sigma ^{2}=\langle x^{2}\rangle -\langle x\rangle ^{2}$, which is
equal to the variance of $x$. This is a measurement of the dispersion (or spread) of the 
particle positions relative to the mean. This also explains why SPH becomes inherently  
unstable when estimating second-order derivatives. For instance, it is well-known that
SPH formulations based on the second-order derivatives of the kernel are highly
sensitive to particle disorder, where dispersion may be enhanced by the presence of
non-uniform velocity fields \cite{Herrera09}. This observation certainly implies that
in practical applications SPH has a limiting second-order convergence.

\section{SPH approximation of first-order derivatives}

The kernel approximation of the gradient of a function can be obtained from Eq. (2)
by replacing $f$ by $\nabla f$ and integrating by parts \cite{Monaghan92}, to yield
\begin{equation}
\langle\nabla f({\bf x})\rangle =\int _{{\cal R}^{3}}f({\bf x}^{\prime})
\nabla W(|{\bf x}-{\bf x}^{\prime}|,h)d{\bf x}^{\prime},
\end{equation}
which can be approximated to $m$th-order accuracy if the following conditions
are satisfied \cite{Liu03}
\begin{equation}
\int _{{\cal R}^{3}}({\bf x}-{\bf x}^{\prime})^{l}\nabla W(|{\bf x}-{\bf x}^{\prime}|,1)
d{\bf x}^{\prime}={\bf I}_{l1},
\end{equation}
where $l=0$,1,2,...$m$ and ${\bf I}_{l1}$ is a tensor of rank $m$ whose components
are all equal to $\delta _{l1}$. For $l=0$, Eq. (22) is equivalent to the condition 
that the kernel function must vanish at the surface of integration, while for $l=1$
it reduces to
\begin{equation}
\int _{{\cal R}^{3}}({\bf x}-{\bf x}^{\prime})\nabla W(|{\bf x}-{\bf x}^{\prime}|,h)
d{\bf x}^{\prime}={\bf I},
\end{equation}
where ${\bf I}$ is the identity matrix. Satisfaction of this condition means that the
isotropy of space is not affected by the SPH kernel approximation \cite{Vignjevic09}
and, as a consequence, angular momentum is preserved \cite{Bonet99}. The same statement
holds for the particle approximation. In other words, the interpolation must be
independent of a rotation of the coordinate axes. In order to see this, let us
consider for simplicity only small rotations so that the coordinates change according
to the transformation
\begin{eqnarray}
{\bf x}^{\prime}&=&{\bf x}-d{\bf w}\times {\bf x}\nonumber\\
&=&{\bf x}-{\bf x}\cdot\nabla (d{\bf w}\times {\bf x}),
\end{eqnarray}
where $d{\bf w}$ is the differential rotation vector. Under solid-body rotation, the
coordinates of a point are independent of the rotation of the coordinate axes and
therefore
\begin{equation}
\langle\nabla (d{\bf w}\times {\bf x})\rangle _{a}=\nabla (d{\bf w}\times {\bf x}),
\end{equation}
or alternatively, using the discretized form of Eq. (21)
\begin{eqnarray}
\langle\nabla (d{\bf w}\times {\bf x})\rangle _{a}&=&\sum _{b=1}^{n}
(d{\bf w}\times {\bf x})_{b}\nabla _{a}W_{ab}\Delta V_{b}\nonumber\\
&=&\sum _{b=1}^{n}[{\bf x}\cdot\nabla (d{\bf w}\times {\bf x})]_{b}\nabla _{a}W_{ab}
\Delta V_{b}\nonumber\\
&=&\nabla (d{\bf w}\times {\bf x})\cdot\sum _{b=1}^{n}{\bf x}_{b}\nabla _{a}W_{ab}
\Delta V_{b},
\end{eqnarray}
which implies that the condition
\begin{equation}
\sum _{b=1}^{n}{\bf x}_{b}\nabla _{a}W_{ab}\Delta V_{b}={\bf I},
\end{equation}
must be satisfied exactly in order to preserve the isotropy of the discretized space
and ensure angular momentum conservation. 

\section{Numerical analysis}

In this section, we describe the results of a series of numerical experiments in two-space
dimensions for two different test functions, namely,
\begin{equation} 
f_{1}(x,y)=\sin\pi x\sin\pi y,
\end{equation}
and
\begin{equation}
f_{2}(x,y)=x^{5/2}\left(20y^{5}+8xy^{3}+x^{2}y^{2}+1\right),
\end{equation}
over the intervals $x\in [0,1]$ and $y\in [0,1]$, in order to check their reproducibility
for both regularly and irregularly distributed particles. We analyze the convergence
rate of the SPH approximations for the functions and their derivatives using five 
different methods: (a) the standard SPH, (b) the CSPM method of Chen et al.
\cite{Chen99a,Chen99b}, (c) the FPM scheme of Liu and Liu \cite{Liu06}, (d) the MSPH
method of Zhang and Batra \cite{Zhang04}, and (e) the methodology recently proposed by
Zhu et al. \cite{Zhu15}, which we label SPH$n$ to distinguish it from the standard SPH.
Similarly, when CSPM and FPM are run with the Wendland function and varying number of 
neighbors we shall label them CSPM$n$ and FPM$n$, respectively.
Unlike the other SPH methods, MSPH
includes second-order derivatives in the Taylor series expansions and therefore it
should restore $C^{2}$ kernel and particle consistency, implying that it should
converge faster than second-order accuracy in contrast to the statement of Eq. (20),
which predicts a second-order limit to the convergence rate of SPH. The cubic
$B$-spline kernel \cite{Monaghan85} with a fixed number of neighbors ($n\approx 13$)
is used for standard SPH, CSPM, FPM, and MSPH, i.e.,  
\begin{equation}
W(q,h)=\frac{15}{7\pi h^{2}}
\left\{
\begin{array}{cll}
\frac{2}{3}-q^{2}+\frac{1}{2}q^{3} & \mbox{if} & 0\leq q<1,\\
\frac{1}{6}(2-q)^{3} & \mbox{if} & 1\leq q<2,\\
0 & \mbox{if} & q\geq 2,
\end{array}
\right.
\end{equation}
where $q=|{\bf x}-{\bf x}^{\prime}|/h$, while a Wendland $C^{4}$ function 
\cite{Wendland95,Dehnen12}
\begin{equation}
W(q,h)=\frac{9}{\pi h^{2}}
\left\{
\begin{array}{cll}
(1-q)^{6}\left(1+6q+\frac{35}{3}q^{2}\right) & \mbox{if} & q<1,\\
0 & \mbox{if} & q\geq 1,
\end{array}
\right.
\end{equation}
is employed for standard SPH, CSPM, and FPM with varying number of neighbors as
the resolution is increased (i.e., schemes SPH$n$, CSPM$n$, and FPM$n$). In these
analyses no boundary treatments are implemented at the borders of the particle
configurations. 

We vary the spatial resolution from 625 to 562500 particles for both distributions
(see first column of Table 1) and measure the convergence rates of the function
estimates in terms of the root mean square error (RMSE)
\begin{equation}
{\rm RMSE}(f)=\sqrt{\frac{1}{N}\sum _{a=1}^{N}\left(f_{a}^{\rm exact}-f_{a}^{\rm num}
\right)^{2}},
\end{equation}
which is closely related to the $L_{2}$-norm. Identical forms to Eq. (32) are also
used to assess the errors in the estimates of the derivatives. Since the SPH errors
are expected to have a normal rather than a uniform distribution, the RMSE will
provide a better representation of the error distribution than other statistical
metrics \cite{Chai14}. Compared to the mean absolute error (MAE), which is more closely
related to the $L_{1}$-norm, the RMSE gives a higher weighting toward large errors in
the sample than the mean average error and therefore it is superior at revealing model
performance differences. In fact, when both metrics are calculated, the RMSE is
always larger than the MAE.
\begin{table}[htbp]
\centering
\caption{Spatial resolution parameters for the calculations}
\begin{tabular}{ccc}
\hline
Number of SPH & Number of & Smoothing length \\
particles & neighbors & \\
$N$ & $n$ & $h$ \\
\hline
~~~625 & ~~213 & 0.342 \\
~~2500 & ~~556 & 0.271 \\
~~5625 & ~~973 & 0.237 \\
~10000 & ~1436 & 0.215 \\
~15625 & ~1933 & 0.200 \\
~22500 & ~2472 & 0.188 \\
~30625 & ~3041 & 0.179 \\
~40000 & ~3648 & 0.170 \\
~62500 & ~4880 & 0.158 \\
~90000 & ~6288 & 0.149 \\
160000 & ~9216 & 0.136 \\
250000 & 12416 & 0.126 \\
562500 & 21328 & 0.110 \\
\hline
\end{tabular}
\end{table}

In the calculations with schemes SPH$n$, CSPM$n$, and FPM$n$, where the Wendland $C^{4}$
kernel is used and the number of neighbors ($n$) is varied with $N$, the quality of the
SPH estimates is analyzed by examining the standard deviation of the SPH evaluated functions
and derivatives as a function of $n$. This provides a measure of the rate at which
the inconsistency of the SPH estimates declines as the number of neighbors increases,
consistently with the expected dependence of the SPH particle discretization error on $n$
\cite{Monaghan92}. For these tests, we use the parameterizations provided by Zhu et 
al. \cite{Zhu15} and allow $h$ to vary with $N$ as $h=N^{-1/6}$. With this choice we
obtain the scaling relations $n\approx 2.81N^{0.675}$ and $h\approx 1.29n^{-0.247}$
for $n$ and $h$, respectively. Figure 1 shows the variation of $h$ with $n$. Based on
a balance between the SPH smoothing and discretization errors, Zhu et al. \cite{Zhu15}
derived the parameterizations $h\propto N^{-1/\beta}$ and $n\propto N^{1-3/\beta}$
for $\beta\in [5,7]$. An intermediate value of $\beta\sim 6$ is appropriate when the
smoothing is performed with the Wendland kernel (31) on a quasi-ordered particle 
configuration. Thus, choosing the proportionality factor of the scaling $h\propto N^{-1/6}$
as exactly unity gives exponents for the dependences of $n$ on $N$ and of $h$ on $n$ that
are slightly larger than 0.5 and $-1/3$, respectively, as suggested by the parameterizations
of Zhu et al. \cite{Zhu15}. However, preliminary tests with exponents close to the
suggested values produced a similar variation of $h$ with $n$ as shown in Fig. 1,
while the convergence rates for SPH$n$, CSPM$n$, and FPM$n$ remained essentially the
same (see next Section). We see that for small values  
of $n$ the smoothing length decreases rapidly as the number of neighbors increases and
then more slowly at large $n$, asymptotically approaching zero as $n\to\infty$ as
needed to restore particle consistency. The second and third columns of Table 1 list 
the number of neighbors and values of the smoothing length, respectively, as they were
obtained from the above scalings.

\subsection{Regularly distributed particles}

As we have discussed before, $C^{m}$ consistency defines the property of a numerical
scheme to reproduce a given field function or distribution to $(m+1)$th-order accuracy.
Therefore, in order to assess the consistency of the several methods we first start by
analyzing their convergence rates on the test functions $f_{1}$ and $f_{2}$ and their
derivatives for a perfectly regular distribution of particles. Table 2 lists the 
calculated convergence rates as obtained from the different SPH methods. For all
methods, the convergence rate corresponds to the slope resulting from a least squares
fitting of the RMSE data points.

\begin{table}[htbp]
\centering
\caption{Convergence rates of RMSE metrics for a regular distribution of SPH particles}
\begin{tabular}{cccccccc}
\hline
Type of function & SPH & CSPM & FPM & MSPH & SPH$n$ & CSPM$n$ & FPM$n$ \\
\hline
$f_{1}$~~~ & ~0~~~ & -0.76 & -1.0~ & -1.76 & -0.97 & -0.92 & -1.0~ \\
$f_{1,x}$~ & +0.51 & -0.77 & -0.77 & -1.0~ & -0.77 & -0.94 & -0.95 \\
$f_{1,y}$~ & +0.51 & -0.77 & -0.77 & -1.0~ & -0.77 & -0.94 & -0.95 \\
$f_{1,xx}$ & ----- & ----- & ----- & -0.75 & ----- & ----- & ----- \\
$f_{1,xy}$ & ----- & ----- & ----- & -1.0~ & ----- & ----- & ----- \\
$f_{1,yy}$ & ----- & ----- & ----- & -0.75 & ----- & ----- & ----- \\
\hline
$f_{2}$~~~ & -0.28 & -0.76 & -1.0~ & -1.76 & -0.78 & -0.92 & -0.98 \\
$f_{2,x}$~ & -0.26 & -0.75 & -0.76 & -1.0~ & -0.77 & -0.91 & -0.93 \\
$f_{2,y}$~ & -0.26 & -0.75 & -0.77 & -1.0~ & -0.77 & -0.91 & -0.93 \\
$f_{2,xx}$ & ----- & ----- & ----- & -0.62 & ----- & ----- & ----- \\
$f_{2,xy}$ & ----- & ----- & ----- & ~0~~~ & ----- & ----- & ----- \\
$f_{2,yy}$ & ----- & ----- & ----- & -0.01 & ----- & ----- & ----- \\
\hline
\end{tabular}
\end{table}

The dependence of the RMSE of the SPH function estimates on the effective number of
particles $N$ is shown in Fig. 2 for $f_{1}$ (left) and $f_{2}$ (right). We see
that standard SPH exhibits very poor convergence rates. For the estimate of $f_{1}$,
the solution converges to $N^{-0.75}$ for the first four resolutions, while at
higher resolutions the RMSE reaches a plateau. Thus, in terms of the absolute error, 
the highest resolution fails to give the best results with standard SPH. Conversely, 
for the estimates of $f_{2}$ the solution converges to $N^{-0.28}$ for all values of
$N$. Faster convergence rates, i.e., $N^{-0.76}$, $N^{-1}$, and $N^{-1.76}$ are obtained
for CSPM, FPM, and MSPH, respectively, for both the estimates of $f_{1}$ and $f_{2}$.
As expected, the fastest convergence is achieved by the MSPH scheme, which is effectively
close to second-order accuracy. Therefore, MSPH is at best $C^{1}$-consistent, while FPM
exhibits $C^{0}$ consistency. When passing from FPM to CSPM the accuracy of SPH degrades
from $N^{-1}$ to $N^{-0.76}$, implying that CSPM is not even restoring full $C^{0}$
consistency. These results validate the arguments presented in Section 3, which show
that SPH has a limiting convergence rate of $N^{-2}$.

In passing, we note that Zhang and Batra \cite{Zhang09} reported convergence rates
for their MSPH scheme of $N^{-3.52}$ when analyzing a sinusoidal test function. Although
they claim to have used an RMSE, we find that using the square of the RMSE, i.e.,
the mean squared error (MSE), which corresponds effectively to an $L_{2}$-norm error,
the convergence rate of MSPH improves to $N^{-3.52}$, which surprinsingly coincides 
with the rate reported by them. This clearly suggests that the convergence rate
found by Zhang and Batra \cite{Zhang09} are more likely based on an MSE than on an
RMSE, and therefore it does not represent a correct estimate of the actual errors.
It is often encountered in the
SPH literature that the MAE, or $L_{1}$-norm error, is also used to determine convergence
rates. However, since the SPH errors are not uniformly distributed, the MAE is likely
to remove information because it gives the same weight to all errors. In contrast,
the RMSE consists of squaring the magnitude of the errors before they are averaged, and
so it gives a relatively higher weight to errors with larger absolute values \cite{Chai14}.
In this sense, the RMSE represents a better statistical metric than the MAE 
because it provides not only a good indicator of the average performance model but also a
better representation of the error distribution.

When standard SPH is run with the Wendland kernel and varied number of neighbors
according to the scaling relation $n\approx 2.81N^{0.675}$,
the convergence rate effectively improves from a near plateau to
$N^{-0.97}$ for the estimate of $f_{1}$ and from $N^{-0.28}$ to $N^{-0.78}$ for
the estimate of $f_{2}$. These rates are comparable to those of CSPM
and FPM with a fixed number ($\sim 13$) of neighbors. Evidently,
increasing the number of neighbors while decreasing the smoothing length with
increasing resolution has the benefit of reducing the SPH discretization errors,
thereby allowing standard SPH to restore $C^{0}$ consistency. However, working
with varying number of neighbors as $N$ is increased has only a moderate effect
on CSPM, whose convergence rate improves from $N^{-0.76}$ to $N^{-0.92}$, and
essentially no effect on FPM, which maintains the same convergence rate as
with fixed $n$. Figure 3 shows the standard deviation of the estimates
of $f_{1}$ and $f_{2}$ as a function of $n$ as obtained with SPH$n$,
CSPM$n$, and FPM$n$. For all three methods, the standard deviation decreases
as we increase the number of neighbors in the SPH estimates. We observe an
approximate $n^{-1}$ trend for all cases, which is appropriate for a quasi-regular
distribution of particles \cite{Zhu15}. This closely follows the expected
dependence of the SPH discretization error on the number of neighbors for a
quasi-ordered pattern \cite{Monaghan92}, $\epsilon\propto\log (n)/n$, which
was then further parameterized by Zhu et al. \cite{Zhu15} as $\epsilon\propto n^{-1}$.
As $n$ is further increased and $h$ is reduced as shown in Fig. 1, the normalization 
condition in discrete form approaches unity. In this limit, the kernel approaches the
Dirac-$\delta$ distribution and the homogeneity of the discretized space is fully
recovered, implying that the field function is exactly reproduced.

The convergence rates of the first-order derivatives ($f_{1,x}$,$f_{1,y}$) and
($f_{2,x}$,$f_{2,y}$) are depicted in Figs. 4 and 5, respectively. The convergence rate 
of standard SPH has a positive slope for the estimates of both the $x$- and $y$-derivatives 
of $f_{1}$, implying that the numerical results diverge from the actual values at all
resolutions. This is a consequence of the RMSE plateauing for the estimate of $f_{1}$
(see Fig. 2). In contrast, the RMSE of the derivative estimates of $f_{2}$
are seen to scale as $N^{-0.26}$ for both derivatives (see Table 2), while for CSPM 
the estimates of the derivatives converge to almost the same rates as the estimate of
the function itself, i.e., as $N^{-0.77}$ for the smoothed derivatives of $f_{1}$
and $N^{-0.75}$ for those of $f_{2}$. As expected, the convergence rate of the derivatives
degrades to less than first-order for the FPM scheme for both test functions and to
first-order for the MSPH scheme, implying $C^{0}$ particle consistency for the first-order
derivatives with this latter method. Almost $C^{0}$ consistency is also restored
for the estimates of the second-order derivatives of $f_{1}$. The estimates of $f_{1,xx}$
and $f_{2,xx}$ share the same convergence rate, i.e., $N^{-0.75}$, while the RMSE
of $f_{1,xy}$ scales as $N^{-1}$. This is not surprising
because $f_{1,xx}=f_{1,yy}=-\pi ^{2}f_{1}$, while $f_{1,xy}=\pi ^{2}\cos\pi x\cos\pi y$.
However, for the polynomial function $f_{2}$, $C^{0}$ consistency is lost as the
RMSE of the estimates of $f_{2,xy}$ and $f_{2,yy}$ follows a plateau for all
resolution runs, while an $N^{-0.62}$ rate is observed for $f_{2,xx}$ (see Table 2).
The results show that SPH$n$, CSPM$n$, and FPM$n$ all exhibit approximately first-order
convergence rates for both the function and its derivatives. This is one advantage
of running with varying number of neighbors over maintaining $n$ fixed with
increasing overall resolution, which restores the same order of consistency for both
the particle estimates of the function and its derivatives.

\subsection{Irregularly distributed particles}

It is well-known that in true SPH simulations, pressure forces mediate between SPH 
particles, which tend to regulate the distribution of neighbors into quasi-regular 
patterns. In some instances, the particle distribution can become highly irregular
as in the case of highly turbulent flows. In either case it is unlikely that
particles will be accommodated in a truly random distribution because the continuity
property of the flow will prevent particles from moving randomly even under highly
non-uniform velocity fields \cite{Monaghan92,Monaghan05}. However, when particles are
arranged in an irregular pattern, a loss of consistency may arise as a consequence
of noise on small scales, resulting in an error which scales as $\sim n^{-1}\ln n$
\cite{Monaghan92}. Moreover, if the degree of particle disorder increases, the
performance of SPH decreases \cite{Herrera13}. Therefore, it is of interest to explore
the behavior of the different SPH methods on particle disorder by maintaining exactly
the same parameters as in the previous section.

Figure 6 depicts the irregular placement of particles for the lowest resolution case (625
particles) on a $[0,1]\times [0,1]$ square. The same degree of disorder was maintained
for all resolutions shown in Table 1. The results of the calculations are condensed in
Table 3, where the convergence rates of the test functions and their derivatives are
shown. Figure 7 displays the dependence of the RMSE for the estimates of $f_{1}$ (left)
and $f_{2}$ (right) on resolution. Comparing the results in Table 3 with those in Table
2 we see that in general the convergence rates of CSPM, FPM, and MSPH slow down for
an irregular pattern. However, schemes SPH$n$, CSPM$n$, and FPM$n$ with varying number
of neighbors are much less sensitive to particle disorder and exhibit almost the same
convergence rates independently of whether the particle distribution is regular or irregular.
The standard deviation of the function estimates against the number of neighbors is displayed 
in Fig. 8. A trend close to $n^{-1}$ is again reproduced for all three methods. Hence,
the small-scale noise induced by particle disorder does not seem to affect the
rate of decay of the standard deviation of the function estimates with $n$. In fact, the
standard deviations show very small discrepancies from an exact partition of unity for
both the regular and irregular point sets. For truly random points, Zhu et al. \cite{Zhu15}
found that the error distribution will follow an $n^{-0.5}$ trend. However, in actual
SPH flow simulations the randomness in the particle distribution is expected to be closer
to an irregular (or quasi-random) sequence, where the spatial particle density is
approximately uniform, rather than to a random one, where large contrasts may exist in
particle density. Therefore, realistic applications may well fall midway between $n^{-1}$
and $n^{-0.5}$ \cite{Zhu15,Herrera13}. However, our results indicate that the dependence
of the standard deviation on $n$ could be more biased toward $n^{-1}$ for irregular
particle configurations.
\begin{table}[htbp]
\centering
\caption{Convergence rates of RMSE metrics for an irregular distribution of SPH particles}
\begin{tabular}{cccccccc}
\hline
Type of function & SPH & CSPM & FPM & MSPH & SPH$n$ & CSPM$n$ & FPM$n$ \\
\hline
$f_{1}$~~~ & -0.07 & -0.59 & -0.99 & -1.53 & -1.02 & -0.92 & -1.0~ \\
$f_{1,x}$~ & +0.52 & -0.62 & -0.62 & -1.10 & -0.78 & -0.94 & -0.94 \\
$f_{1,y}$~ & +0.51 & -0.62 & -0.62 & -1.13 & -0.78 & -0.94 & -0.94 \\
$f_{1,xx}$ & ----- & ----- & ----- & -0.59 & ----- & ----- & ----- \\
$f_{1,xy}$ & ----- & ----- & ----- & -0.65 & ----- & ----- & ----- \\
$f_{1,yy}$ & ----- & ----- & ----- & -0.72 & ----- & ----- & ----- \\
\hline
$f_{2}$~~~ & -0.18 & -0.59 & -0.99 & -1.53 & -0.77 & -0.93 & -1.0~ \\
$f_{2,x}$~ & -0.13 & -0.67 & -0.76 & -1.08 & -0.78 & -0.93 & -0.94 \\
$f_{2,y}$~ & -0.12 & -0.65 & -0.64 & -1.11 & -0.78 & -0.92 & -0.94 \\
$f_{2,xx}$ & ----- & ----- & ----- & -0.73 & ----- & ----- & ----- \\
$f_{2,xy}$ & ----- & ----- & ----- & -0.02 & ----- & ----- & ----- \\
$f_{2,yy}$ & ----- & ----- & ----- & -0.03 & ----- & ----- & ----- \\
\hline
\end{tabular}
\end{table}

For completeness, Figs. 9 and 10 depict the convergence rates of the smoothed derivatives 
of $f_{1}$ and $f_{2}$, respectively. The trends are similar to those displayed in
Figs. 4 and 5 for a regular distribution of particles, except for the runs with MSPH
for which the variations of the RMSE with $N$ shows up and down turns, which are
considerably more pronounced for the estimates of $f_{1}$ than for those of $f_{2}$.
This zigzag behavior is indicative of MSPH being more sensitive to particle disorder
when estimating derivatives. Similar erratic behaviors are also observed for the RMSE of
the estimates of the second-order derivatives, except for the solutions of $f_{2,xy}$ and 
$f_{2,yy}$, which fail to converge because their RMSE follows a plateau for all resolution
runs. As noted earlier, the convergence rates of the smoothed estimates of the 
second-order derivatives given in Tables 2 and 3 correspond to least squares fittings
of the numerical errors. 

In general, the convergence rates of the estimates of the first- and second-order
derivatives slow down when the smoothing is performed on an irregular distribution of
points. This is the case of standard SPH, CSPM, FPM, and MSPH with a
fixed number of neighbors. However, if we compare the last three columns of Tables
2 and 3, we see that there are little differences between the convergence rates of
SPH$n$, CSPM$n$, and FPM$n$ for regular and irregular particle distributions.     
As $n$ is increased with resolution, the discretization error, which scales
as $n^{-1}$, decreases, making the calculations with varied number of neighbors much
less sensitive to particle disorder. Similarly to the regular configuration, the
convergence rate of standard SPH for the estimates of $f_{1}$ improves from $N^{-0.07}$
to first-order accuracy ($N^{-1.02}$), while for $f_{2}$ the speed of convergence
increases from $N^{-0.18}$ to $N^{-0.77}$, implying that $C^{0}$ consistency of
standard SPH is restored when working with varied number of neighbors. First-order
convergence rates are also observed with SPH$n$ for the estimates of the first-order
derivatives. Full $C^{0}$ consistency is also restored with CSPM$n$ for both the 
estimates of $f_{1}$ and $f_{2}$ and their first derivatives. Thus, the same order
of consistency is maintained for the approximation of the function and its derivatives
even for a disordered point set. This is also seen when comparing FPM with FPM$n$,
where the accuracy of the approximation for the derivatives becomes close to
first order for both test functions for the latter scheme. Hence, as long as $N\to\infty$,
$h\to 0$, and $n\to\infty$ we expect that the SPH estimates of the derivatives 
converge essentially at the same rate as the estimate of the function. In this limit,
the approximation becomes insensitive to particle disorder.

\section{Conclusions}

In this paper, we have re-examined the problem of kernel and particle consistency
in the smoothed particle hydrodynamics (SPH) method. In particular, we first
demonstrate with a simple observation that any kernel function that is suitable
for SPH interpolation can be scaled with respect to the smoothing length in such a
way that the kernel normalization condition as well as the family of consistency
relations become independent of the smoothing length. This has important implications
on the issue of particle consistency in that the discrete summation form of the
integral consistency relations will only depend on the number of neighboring particles
within the kernel support. While this result was previously derived from detailed
SPH error analyses by Vaughan et al. \cite{Vaughan08} and Read et al. \cite{Read10},
we find as a further implication of our analysis that $C^{2}$ kernel consistency is
difficult to achieve in actual SPH simulations due to an intrinsic diffusion, which
is closely related to the inherent dispersion of the particle positions relative to
the mean. This implies that in practical applications SPH has a limiting second-order
convergence rate. 

Numerical experiments with suitably chosen test functions in two-space dimensions
validate our findings for a number of well-known SPH methods, namely the standard
SPH, the CSPM method of Chen et al. \cite{Chen99a,Chen99b}, the FPM scheme of Liu
and Liu \cite{Liu06}, the MSPH method of Zhang and Batra \cite{Zhang04}, and the
recently proposed methodology of Zhu et al. \cite{Zhu15}, where no corrections
are required and full consistency is restored by allowing the number of neighbors
to increase and the smoothing length to decrease with increasing spatial resolution.
In particular, we find that when using the root mean square
error (RMSE) as a model evaluation statistics, CSPM and FPM converge to only 
first-order accuracy, while MSPH, which was previously thought to converge to 
better than third order, is actually close to second order. This result implies
that both CSPM and FPM have at best $C^{0}$ consistency, while MSPH has $C^{1}$
consistency. When standard SPH is run with varied number of neighbors according
to the method of Zhu et al. \cite{Zhu15}, $C^{0}$ consistency is fully restored for
both the estimates of the function and its derivatives. The same is observed for
the CSPM and FPM methods. The results also show that with varying number of
neighbors the order of consistency is not affected by particle disorder in
contrast to the case where the number of neighbors is kept fixed. Although the
method of Zhu et al. \cite{Zhu15} restores only $C^{0}$ consistency for the
spatial resolutions that are attainable with the use of present-day computers,
it has the advantage of being insensitive to particle disorder and yielding
estimates of the function and its derivatives that converge essentially
at the same rate. However, in terms of the computational cost there remains the 
question on the feasibility of the method for practical applications, where 
restoring full $C^{0}$ consistency in highly resolved calculations will demand
using quite a large number of neighbors compared to conventional SPH methods
where typical choices of $n$ lie in the range $\sim 12$--27 in two-space dimensions
and $\sim 33$--64 in three-space dimensions.

The results of the present analysis not only highlight the complexity of error
behavior in SPH, but also show that restoring $C^{2}$ consistency, or
equivalently, achieving an accuracy higher than second order, still remains
a challenge.

\section*{Acknowledgement} One of us, F. P.-P. is grateful to ABACUS for financial
support during his visit to the Department of Mathematics of Cinvestav, Mexico. 
This work is partially supported by ABACUS, CONACyT grant EDOMEX-2011-C01-165873
and by the Departamento de Ciencias B\'asicas e Ingenier\'{\i}a (CBI) of the
Universidad Aut\'onoma Metropolitana--Azcapotzalco (UAM-A) through internal
funds.

\section*{References}

\bibliography{mybibfile}

\clearpage

\begin{figure}
\centering
\includegraphics[width=10cm]{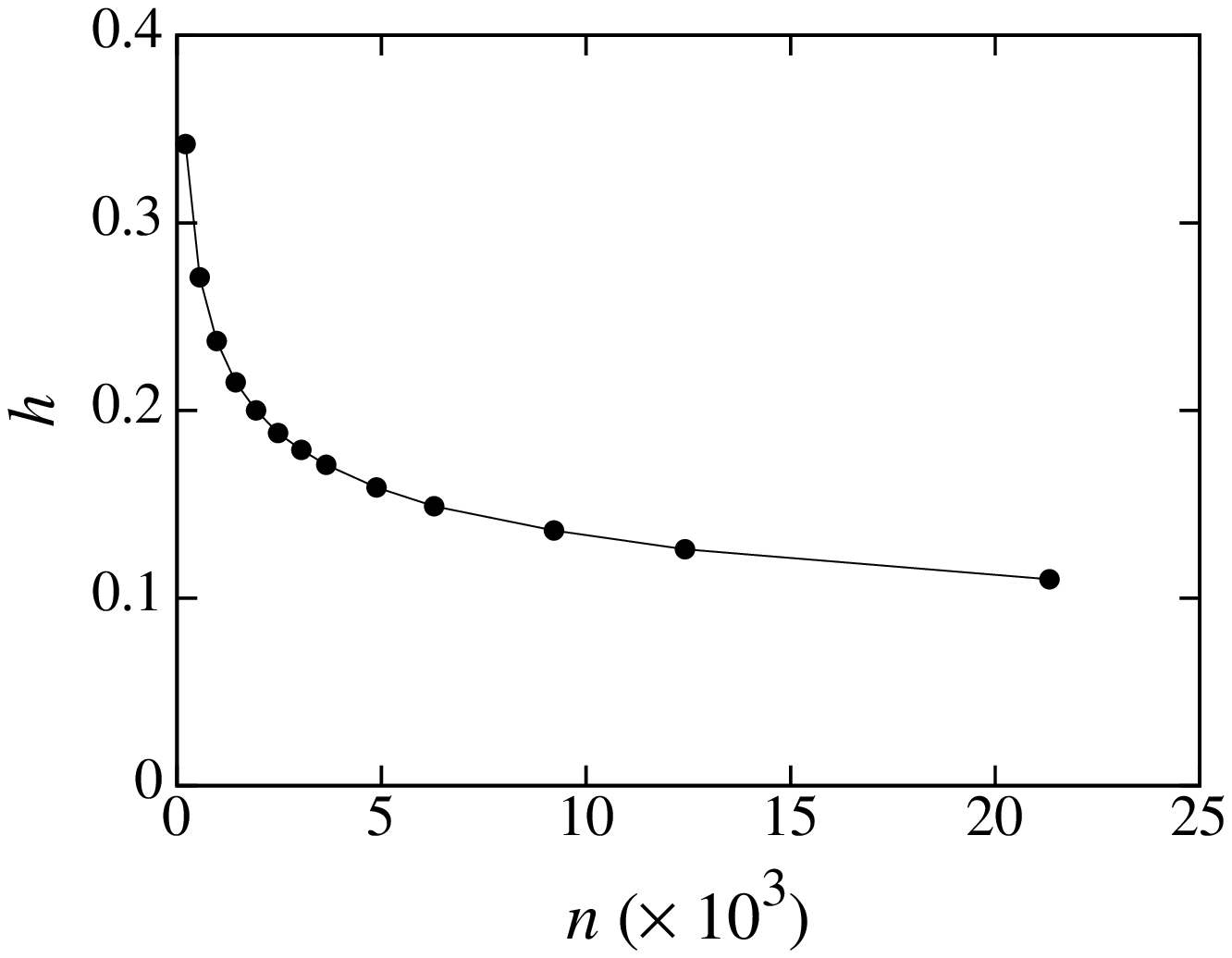}
\caption{Dependence of the smoothing length, $h$, on the number of neighbors, $n$,
as calculated using the scaling relation $h\approx 1.29n^{-0.247}$.} 
\end{figure}
\begin{figure}
\centering
\includegraphics[width=10cm]{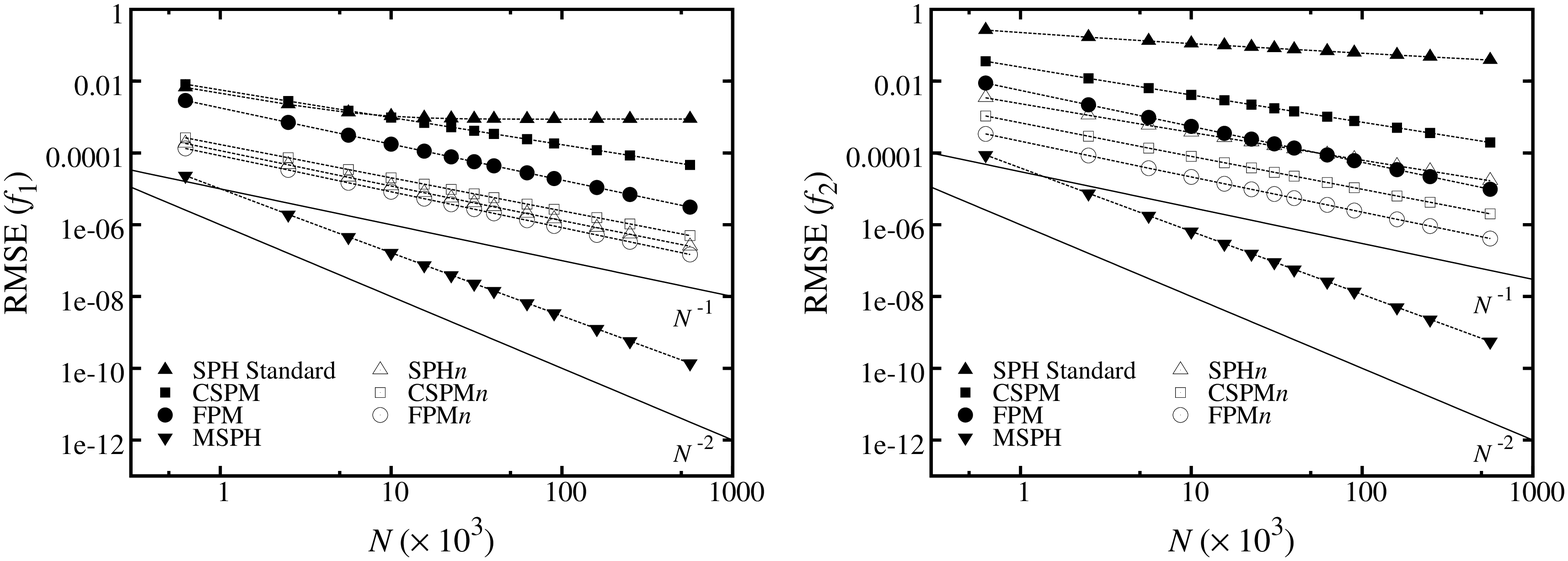}
\caption{Convergence rates of the different SPH methods for the estimates of $f_{1}$ 
(left) and $f_{2}$ (right) as a function of the total number of particles, $N$, as
obtained for a perfectly regular distribution of particles. The $N^{-1}$ and $N^{-2}$ 
trends are shown for comparison. The slopes of the straight lines for each method
are given in Table 2.}
\end{figure}
\begin{figure}
\centering
\includegraphics[width=10cm]{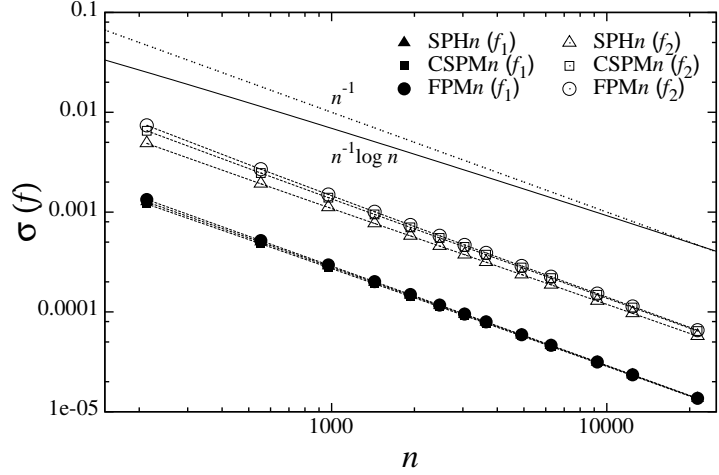}
\caption{Standard deviation of the estimates of $f_{1}$ and $f_{2}$ as a
function of $n$ as obtained for a perfectly regular configuration of particles.
The $n^{-1}$ and $n^{-1}\log n$ trends are displayed for comparison. The slopes of the
straight lines are: $-0.962$ (SPH$n$), $-1.0$ (CSPM$n$), $-1.02$ (FPM$n$) for the 
estimates of $f_{1}$ (filled markers) and $-0.985$ (SPH$n$), $-0.976$ (CSPM$n$), $-0.995$
(FPM$n$) for the estimates of $f_{2}$ (empty markers).}
\end{figure}
\begin{figure}
\centering
\includegraphics[width=10cm]{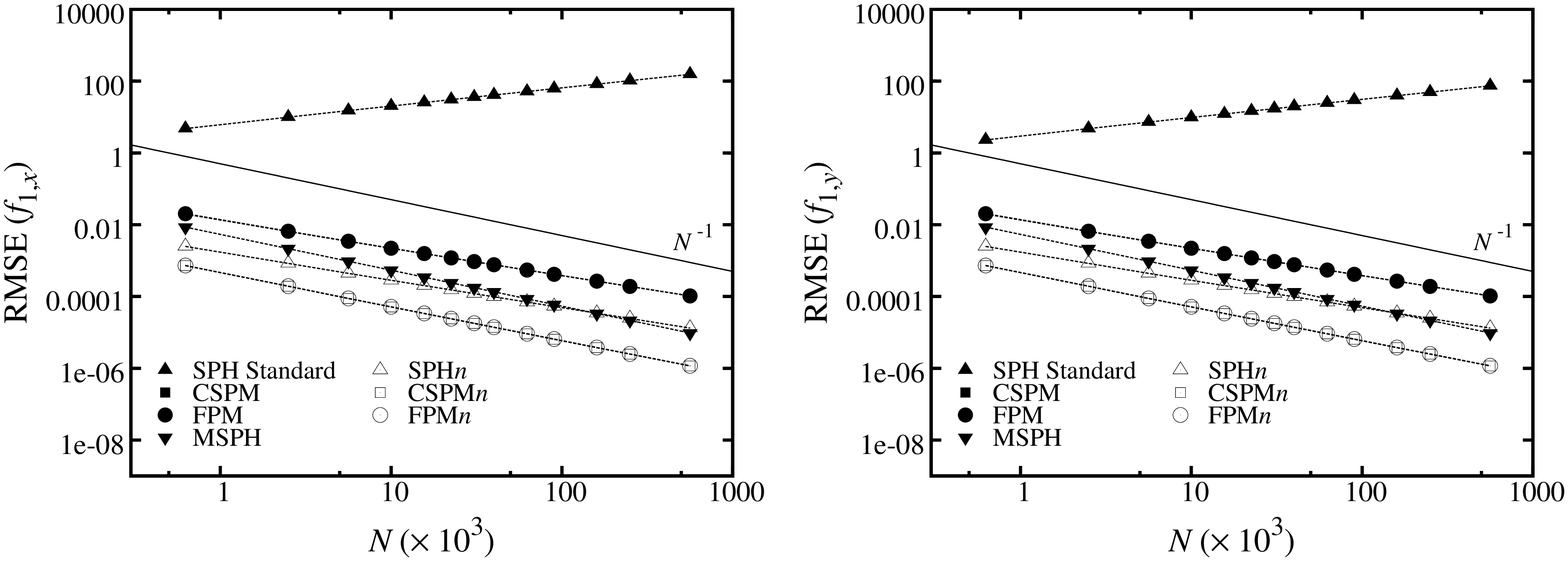}
\caption{Convergence rates of the different SPH methods for the estimates of the 
first-order derivatives of $f_{1}$, $f_{1,x}$ (left) and $f_{1,y}$ (right), as a
function of the total number of particles, $N$, for a perfectly regular distribution
of particles. The $N^{-1}$ trend is shown for comparison. The slopes of the straight
lines are listed in Table 2.}  
\end{figure}
\begin{figure}
\centering
\includegraphics[width=10cm]{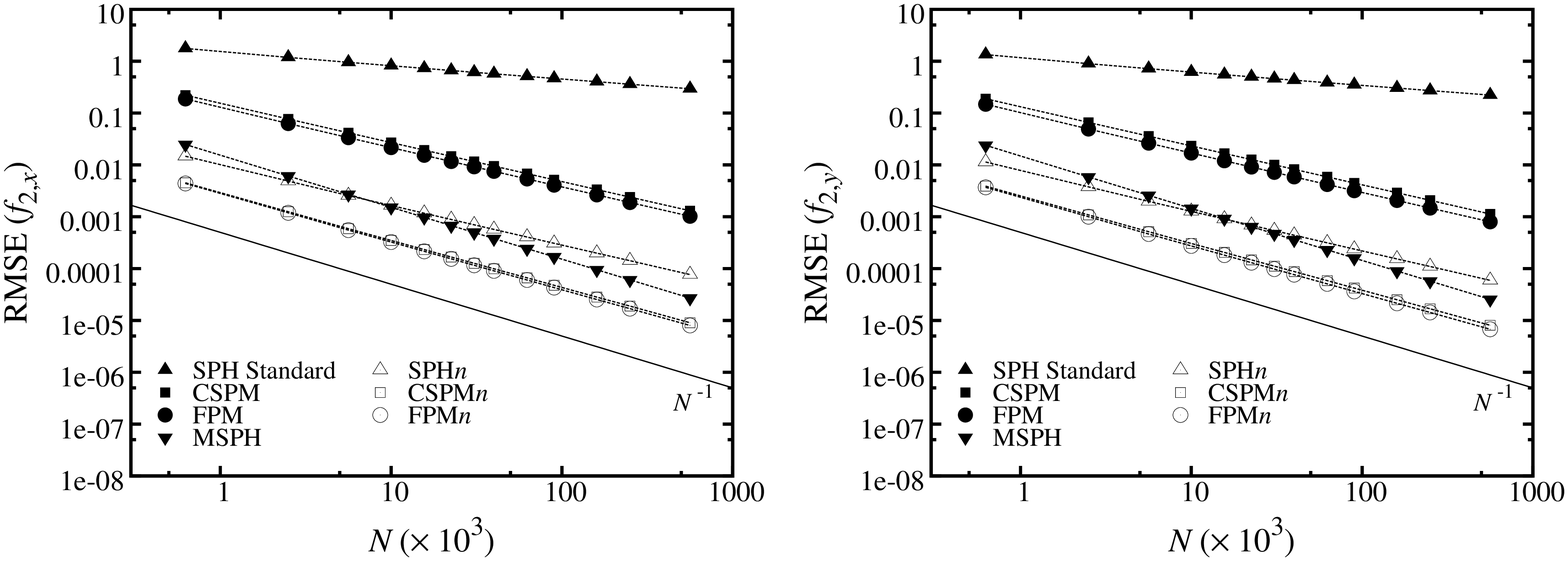}
\caption{Convergence rates of the different SPH methods for the estimates of the 
first-order derivatives of $f_{2}$, $f_{2,x}$ (left) and $f_{2,y}$ (right), as a
function of the total number of particles, $N$, for a perfectly regular distribution
of particles. The $N^{-1}$ trend is shown for comparison. The slopes of the straight
lines are listed in Table 2.}  
\end{figure}
\begin{figure}
\centering
\includegraphics[width=10cm]{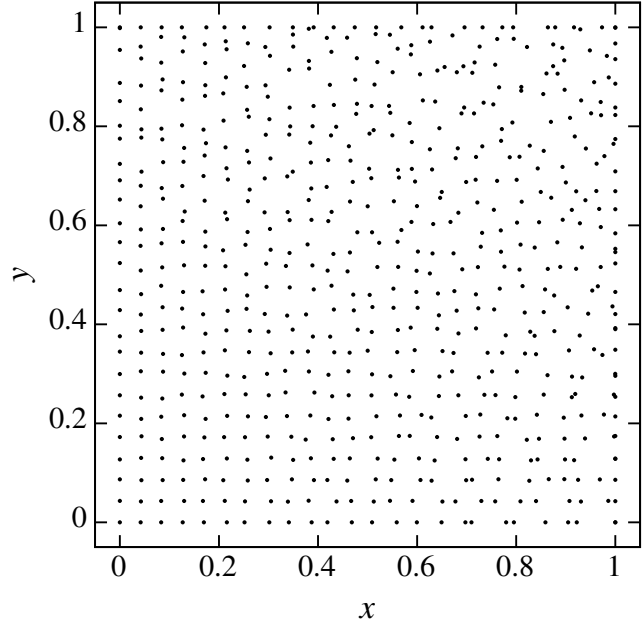}
\caption{Irregular particle distribution for the lowest resolution case ($N=625$)
in a two-dimensional domain $[0,1]\times [0,1]$. The same degree of particle
disorder is maintained for all higher resolutions.}
\end{figure}
\begin{figure}
\centering
\includegraphics[width=10cm]{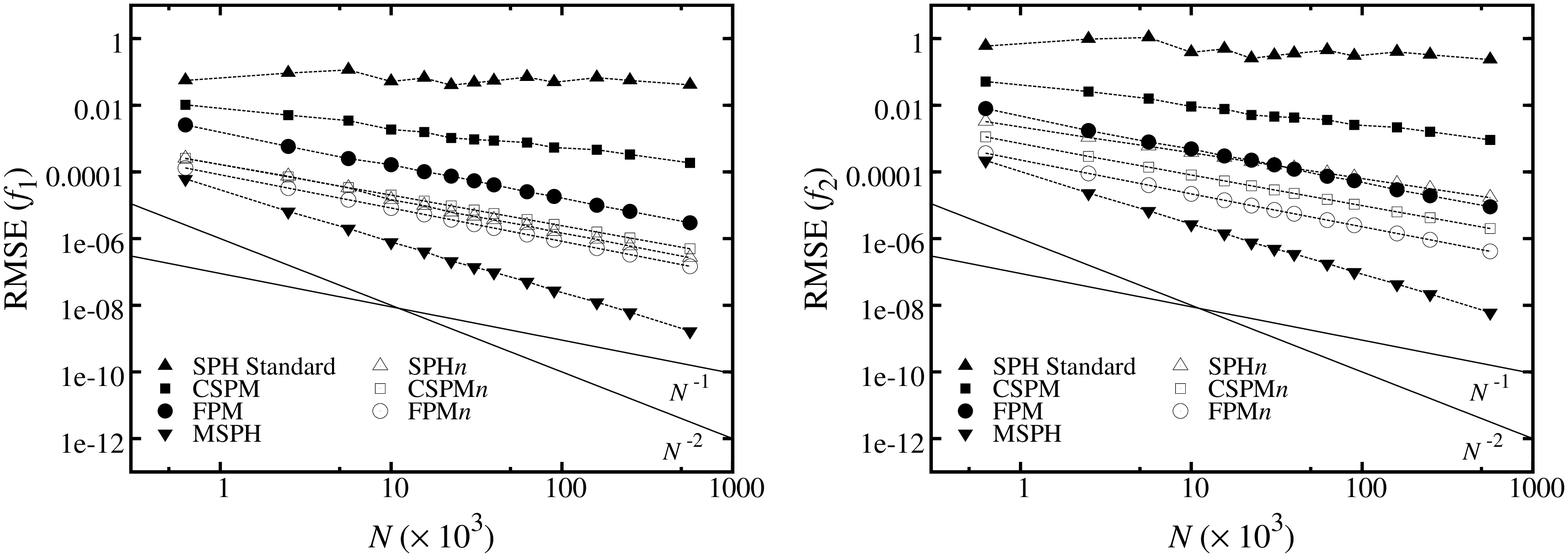}
\caption{Convergence rates of the different SPH methods for the estimates of $f_{1}$ 
(left) and $f_{2}$ (right) as a function of the total number of particles, $N$, for an
irregular distribution of particles. The $N^{-1}$ and $N^{-2}$ trends are shown for
comparison. The slopes of the straight lines for each method are given in Table 3.}
\end{figure}
\begin{figure}
\centering
\includegraphics[width=10cm]{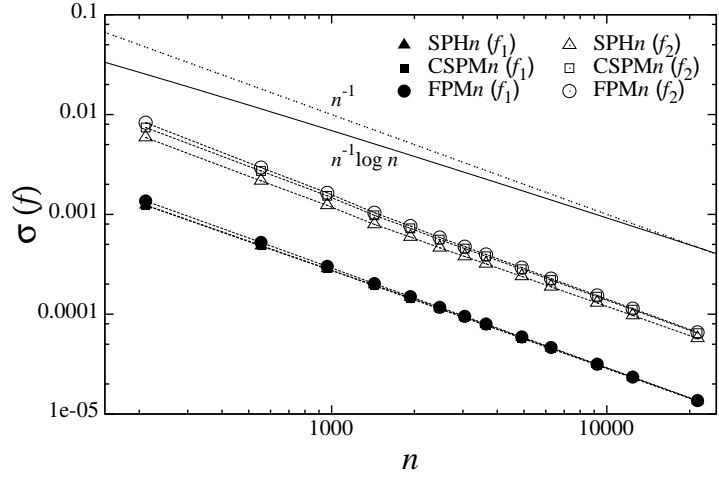}
\caption{Standard deviation of the estimates of $f_{1}$ and $f_{2}$ as a
function of $n$ as obtained for an irregular configuration of particles.
The $n^{-1}$ and $n^{-1}\log n$ trends are displayed for comparison. The slopes of the
straight lines are: $-0.978$ (SPH$n$), $-0.982$ (CSPM$n$), $-0.999$ (FPM$n$) for the
estimates of $f_{1}$ (filled markers) and $-1.0$ (SPH$n$), $-1.03$ (CSPM$n$), $-1.05$
(FPM$n$) for the estimates of $f_{2}$ (empty markers).} 
\end{figure}
\begin{figure}
\centering
\includegraphics[width=10cm]{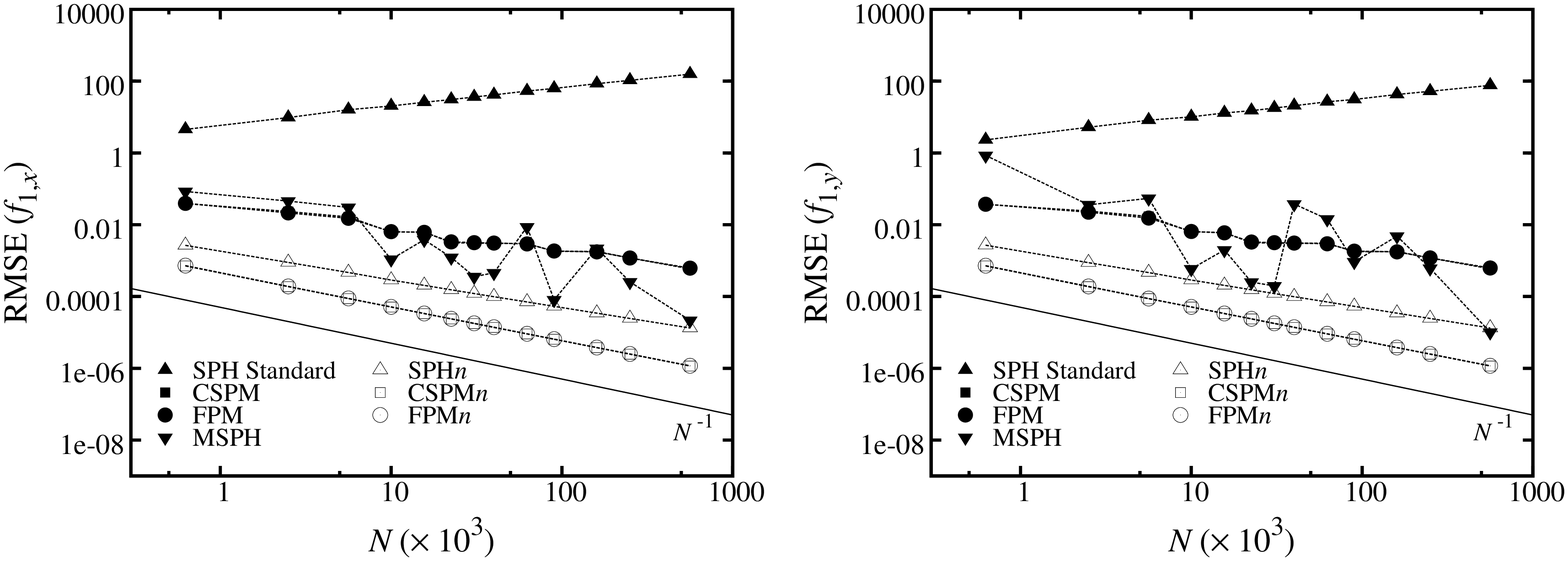}
\caption{Convergence rates of the different SPH methods for the estimates of the 
first-order derivatives of $f_{1}$, $f_{1,x}$ (left) and $f_{1,y}$ (right), as a
function of the total number of particles, $N$, for an irregular distribution of
particles. The $N^{-1}$ trend is shown for comparison. The slopes of the straight
lines are given in Table 3.}
\end{figure}
\begin{figure}
\centering
\includegraphics[width=10cm]{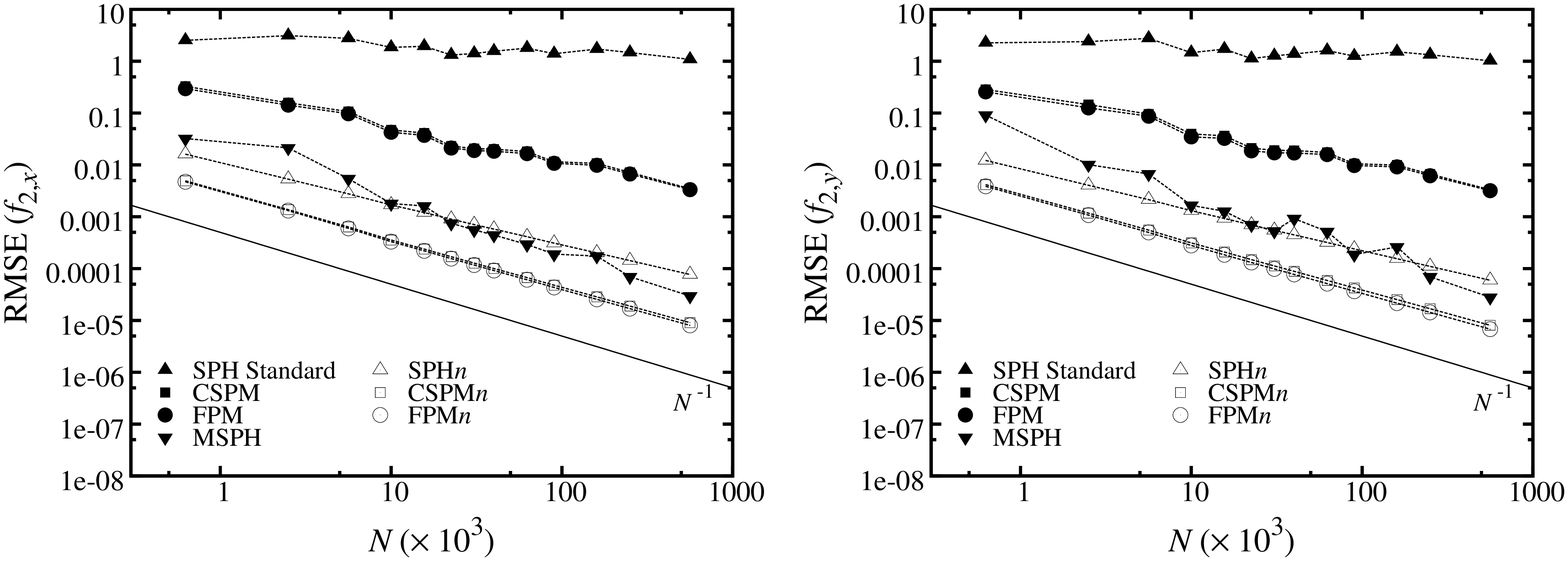}
\caption{Convergence rates of the different SPH methods for the estimates of the 
first-order derivatives of $f_{2}$, $f_{2,x}$ (left) and $f_{2,y}$ (right), as a
function of the total number of particles, $N$, for an irregular distribution of
particles. The $N^{-1}$ trend is shown for comparison. The slopes of the straight
lines are given in Table 3.}
\end{figure}

\end{document}